\documentclass[10pt]{article}

\usepackage{amsmath}
\usepackage{amssymb}
\usepackage{amsthm}
\usepackage{graphicx}

\newtheorem{theorem}{Theorem}

\newtheorem{defin}{Definition}[section]
          
          \newtheorem{teo}{Theorem}[section]
          
          \newtheorem{con}{Conjecture}
          \newtheorem{cond}{Condition}
          \newtheorem{prop}[teo]{Proposition}
          \newtheorem{lem}{Lemma}[section]
          \newtheorem{rmk}[teo]{Remark}
          \newtheorem{cor}{Corollary}[section]
          
          \newcommand{\beq}{\begin{equation}}
          \newcommand{\eeq}{\end{equation}}
          \newcommand{\beqn}{\begin{eqnarray}}
          \newcommand{\beqnn}{\begin{eqnarray*}}
          \newcommand{\eeqn}{\end{eqnarray}}
          \newcommand{\eeqnn}{\end{eqnarray*}}
          \newcommand{\bprop}{\begin{prop}}
          \newcommand{\eprop}{\end{prop}}
          \newcommand{\bteo}{\begin{teo}}
          \newcommand{\bcor}{\begin{cor}}
          \newcommand{\ecor}{\end{cor}}
          \newcommand{\bcon}{\begin{con}}
          \newcommand{\econ}{\end{con}}
          \newcommand{\bcond}{\begin{cond}}
          \newcommand{\econd}{\end{cond}}
          \newcommand{\eteo}{\end{teo}}
          \newcommand{\brm}{\begin{rmk}}
          \newcommand{\erm}{\end{rmk}}
          \newcommand{\blem}{\begin{lem}}
          \newcommand{\elem}{\end{lem}}
          \newcommand{\ben}{\begin{enumerate}}
          \newcommand{\een}{\end{enumerate}}
          \newcommand{\bei}{\begin{itemize}}
          \newcommand{\eei}{\end{itemize}}
          \newcommand{\bdf}{\begin{defin}}
          \newcommand{\edf}{\end{defin}}

\newcommand{\su}{{\tilde  S^{\tau}}}
\newcommand{\sd}{{\tilde S^{\tau'}}}
\newcommand{\sud}{{\tilde S^{\tau}_{d}}}
\newcommand{\sdd}{{\tilde S^{\tau'}_{d}}}
\newcommand{\Ss}{{\tilde S_{s}}}

          \newcommand{\nn}{\nonumber}

          \newcommand{\fr}{\frac}
          \renewcommand{\r}{{\mathbb R}}
          \newcommand{\br}{\bar{\mathbb R}}
          \newcommand{\Z}{{\mathbb Z}}

          \newcommand{\R}{{\mathbb R}}
          
          \newcommand{\E}{{\mathbb E}}
          
          \renewcommand{\P}{{\mathbb P}}
          \newcommand{\N}{{\mathbb N}}
          
          \renewcommand{\S}{{\cal S}}

          \newcommand{\W}{{\cal W}}

          \newcommand{\h}{{\cal H}}
          \newcommand{\f}{{\cal F}}

          \renewcommand{\a}{\alpha}
          \renewcommand{\b}{\beta}
          \newcommand{\g}{\gamma}

          \renewcommand{\d}{\delta}
          \newcommand{\D}{\Delta}
          \newcommand{\e}{\epsilon}

          \newcommand{\s}{\sigma}
          \renewcommand{\o}{\Pi}

\newcommand{\txt}{\textrm}

\newcommand{\btt}{\begin{theorem}}
\newcommand{\ett}{\end{theorem}}

\newcommand{\daw}{\downarrow}
\newcommand{\uaw}{\uparrow}
\newcommand{\raw}{\rightarrow}

\newcommand{\be}{\begin{equation}}
\newcommand{\ee}{\end{equation}}

          \newcommand\sqr{\vcenter{
          \hrule height.1mm
          \hbox{\vrule width.1mm height2.2mm\kern2.18mm\vrule width.1mm}
          \hrule height.1mm}}        
\begin{document}

\title{Exceptional Times for the Dynamical Discrete Web}

\author{
{\bf L.~R.~G.~Fontes}
{\small Instituto de Matem\'{a}tica e Estat\'{i}stica,
Universidade de S\~{a}o Paulo, Brazil}\\
\and
{\bf C.~M.~Newman}
{\small Courant Inst.~of Mathematical Sciences, NYU, New York, NY 10012}
\and
{\bf K.~Ravishankar}
{\small  Dept.~of Mathematics, SUNY College at New Paltz, New Paltz, NY 12561}
\and
{\bf E.~Schertzer}
{\small  Courant Inst.~of Mathematical Sciences, NYU, New York, NY 10012}
}

\date{}
\maketitle

\begin{abstract}
The dynamical discrete web (DyDW),
introduced in recent work of Howitt and Warren,
is a system of coalescing simple symmetric
one-dimensional random walks which evolve in an extra continuous dynamical
time parameter $\tau$. The evolution is by independent updating of the underlying
Bernoulli variables indexed by discrete space-time that define the discrete
web at any fixed $\tau$. In this paper, we study the existence 
of exceptional (random) values of $\tau$ where the paths of the web
do not behave like usual random walks and the Hausdorff dimension
of the set of exceptional such $\tau$. 
Our results are motivated by
those about exceptional times for dynamical percolation
in high dimension by H\"{a}ggstrom, Peres and Steif, 
and in dimension two by Schramm and Steif. The exceptional behavior of
the walks in the DyDW 
is rather different from the situation  
for the dynamical random walks of Benjamini,
H\"{a}ggstrom, Peres and Steif.
For example, we prove that 
the walk from the origin $S^\tau_0$ violates the law of the iterated logarithm (LIL)
on a set of 
$\tau$ of Hausdorff dimension one.
We also discuss how these and
other results extend to the dynamical Brownian web,
the natural scaling limit of the DyDW.
\end{abstract}

\newpage
\section{Introduction}
\label{intro}

In this paper, we present a number of results concerning a dynamical
version of coalescing random walks,
which was recently introduced in~\cite{HW07}.
Our results concern sets of 
dynamical times of Hausdorff dimension less than or equal
to one (and of zero Lebesgue measure) where
the system of coalescing walks behaves exceptionally. The results are analogous to
and were motivated by the model of dynamical percolation and
its exceptional times~\cite{HPS97,SS05}. In this section, we define
the basic model treated in this paper, which we call 
the dynamical discrete web (DyDW),
recall some facts about dynamical percolation, and then briefly describe
our main results. The justification for calling this model a discrete web
is that there is a natural scaling limit, the dynamical Brownian web (DyBW), 
which 
was proposed by Howitt and Warren in~\cite{HW07}
and completely constructed 
in~\cite{NRS08}. As we shall explain (see Section \ref{last-section2}),
the exceptional times results for the DyDW 
extend to the continuum DyBW.

We 
note that exceptional times for other dynamical
versions of random walks in various spatial dimensions
have been studied in \cite{BHPS03,Hoff05,AH06} and elsewhere, but,
as we shall see, 
these are quite different from the dynamical random walks of the DyDW.

\bigskip

{\bf The Discrete Web} 

The discrete web is a collection of coalescing 
one-dimensional simple random walks starting 
from every
point in the discrete 
space-time domain 
$\Z_{even}^2=\{(x,t)\in\Z^2: x+t \ \txt{is even}\}$. 
The Bernouilli percolation-like structure is highlighted by defining $\xi_{x,t}$ for
$(x,t)\in \mathbb{Z}_{even}^2$ to be the increment of the random walk
at location $x$ at time $t$. These Bernoulli variables
are symmetric and independent and the paths of all the coalescing 
random
walks can be reconstructed by assigning to each point $(x,t)$
an arrow from  $(x,t)$ to $\{x+\xi_{x,t},t+1\}$
and
considering
all the paths starting from arbitrary points in $\mathbb{Z}_{even}^2$
that follow the arrow configuration $\aleph$.

\bigskip

{\bf The Dynamical Discrete Web}

In the DyDW, there is, in addition
to the random
walk discrete time parameter, 
an additional (continuous)
dynamical time
parameter $\tau$. The system 
starts at $\tau=0$ as an ordinary
DW and then evolves in $\tau$ 
by randomly switching the direction of each
arrow at a fixed rate 
independently of all other arrows. 
We will generally do the switching by having at each 
$(x,t)\in\Z_{even}^2$ a Poisson
clock ring at rate one and then reset 
the direction of the arrow at random; thus 
the rate of switching will be
$1/2$.
This amounts to 
extending the percolation substructure $\xi_z^0$
to time varying functions $\xi_z^\tau$ defining
a (right continuous) dynamical arrow configuration 
$\tau\rightsquigarrow\aleph(\tau)$ and 
$W(\tau)$,
the dynamical discrete web at time $\tau$, is defined as the web
constructed from  $\aleph(\tau)$. 

If one follows the arrows
starting from the (space-time) origin $(0,0)$, the dynamical path 
$S_0^{\tau}$ begins
at $\tau=0$ as a simple symmetric random walk and then evolves dynamically
in $\tau$. At any fixed time $\tau$,
$S^\tau_0$ has
the same law as at time $\tau=0$.
As a consequence, if $\mu$ is the probability distribution of a
simple symmetric random walk starting from the origin and $A$ is any 
event
with $\mu(A)=1$,  we have
for any deterministic $\tau$ that $\P(S^{\tau}_0\in A)=1$. By 
a straightforward application of 
Fubini's Theorem this implies that
\begin{equation}
\label{rrr}
\P(S^{\tau}_0\in A \ \txt{for Lebesgue a.e. $\tau$})=1.
\end{equation}
Following \cite{BHPS03}, for any event such that (\ref{rrr})
holds, a natural question is whether (\ref{rrr}) can be strengthened to
\begin{equation}
 \P(S^{\tau}\in A \ \txt{for all $\tau\geq 0$})=1,
\end{equation}
i.e., {\it do there exist some exceptional times $\tau$
at which $S_0^\tau$ violates
some almost sure 
properties of the standard random walk?} or stated differently, 
is the  random walk sensitive to the dynamics introduced on the DW?

\bigskip

{\bf Analogies With Dynamical Percolation}

Similar questions have been investigated in percolation.
Static (site) 
percolation models are defined also in terms of independent Bernoulli
variables 
$\xi_z^0$, indexed by points $z$ in some $d$-dimensional lattice,
which in general are asymmetric with parameter $p$. There is a critical value
$p_c$ when the system has a transition from having an infinite cluster (connected
component) with probability zero to having one with probability one. It is
expected that at $p=p_c$ there are no infinite clusters and this is proved for
$d=2$ and for high $d$ (see, e.g., \cite{G89}).
In dynamical percolation, one
extends $\xi_z^0$ to time varying functions $\xi_z^\tau$, as in the case
of coalescing walks, except that the transition rates for the jump processes
$\xi_z^\tau$ are chosen to have the critical asymmetric $(p_c,1-p_c)$
distribution to be invariant. The question raised 
in~\cite{HPS97} 
was whether there were exceptional times when an infinite cluster (say,
one containing the origin) occurs, 
even though this does not occur at deterministic
times. This was answered negatively in~\cite{HPS97} for large $d$
and, more remarkably, was answered positively
by Schramm and Steiff for $d=2$ in~\cite{SS05},
where they further obtained upper and lower bounds on the Hausdorff
dimension (as a subset of the dynamical time axis) of these exceptional times.
In \cite{GaPS08}, the exact Hausdorff dimension was obtained. 

\bigskip

{\bf Main Results}

We apply in this paper the approaches used for dynamical percolation
to the dynamical discrete web. Although we restrict attention to
one-dimensional random walks
whose paths are in two-dimensional
space-time and hence analogous to $d=2$ dynamical percolation,
we use
both the high~$d$ and $d=2$ methods
of~\cite{HPS97,SS05}.

{\bf --Tameness.} A natural initial question 
is whether there might be exceptional
dynamical times~$\tau$ for which 
the walk from the origin $S^\tau_0(t)$ is transient
(say to $+\infty$). Our first main result (see Theorem~\ref{th1}),
modeled after the
high-$d$ dynamical percolation results of~\cite{HPS97}, is that
there are no such exceptional times.

{\bf --Existence of Exceptional Times.} 
For a simple symmetric
random walk $S$, it is well known that $\liminf_{t\uaw\infty} S(t)/\sqrt{t}=-\infty$ a.s. (and, of course, $\limsup_{t\uaw\infty} S(t)/\sqrt{t}=+\infty$ a.s.). In the following, we will 
say that a path is subdiffusive if it violates this a.s. property of the
standard random walk.

\bdf\label{K-subd}{\bf [$K^{+}$-subdiffusivity]}
Let $K\in(0,\infty)$.
A path $\pi$ starting at $x=0$ at time $t=0$ is said to be $K^+$-subdiffusive iff there exists $j\geq0$
such that
\begin{equation}
 \forall t>0, \ \ \pi(t)\geq -j-K\sqrt{t}.
\end{equation}
We say that $\pi$ is subdiffusive iff there exists $K\in(0,\infty)$ such that $\liminf_{t\uaw\infty}\pi(t)/\sqrt{t}\geq-K$ or $\limsup_{t\uaw\infty}\pi(t)/\sqrt{t}\leq K$.
\edf
In Proposition~\ref{violation}, we prove that for $K$ large enough, there is a
strictly  positive 
probability for having a dynamical
time $\tau\in[0,1]$ at which $S^\tau_0$ is
$K^+$-subdiffusive.
Propositions~\ref{lower} and ~\ref{upper} gives lower and upper bounds
on the (deterministic) Hausdorff dimension of these exceptional times 
in $[0,\infty)$.
Interestingly, the Hausdorff dimensions depend nontrivially on the constant $K$ so that the
dimension tends to zero (respectively, one) as $K\to0$ (respectively,
$K\to \infty$). In particular, as a direct consequence of Proposition \ref{lower}, 
we obtain the following result.
\bteo
The set of times 
$\tau \in [0,\infty)$ at which
$S_0^\tau$ is subdiffusive has Hausdorff dimension one. 
Hence, the set of exceptional times for
the  
law of the iterated logarithm (LIL) also  has Hausdorff dimension one.
\eteo
Since a set of exceptional times has zero Lebesgue measure (see (\ref{rrr})), we see 
that the set
of exceptional times for the LIL (or for subdiffusivity) is in a sense 
as large as it can be. This is strikingly in contrast with the dynamical one-dimensional
random walk of
\cite{BHPS03} where there are no exceptional times for which
the LIL fails (in \cite{BHPS03}, we recall that the analogue of $S^\tau_0$
is simply defined as
\begin{equation}
\label{other-dy}
\bar S^{\tau}_0(n)=\sum_{i=1}^{n} X_i^\tau,
\end{equation}
where $\{X_i^\tau\}_i$ are 
independent $\{-1,+1\}$-valued Markov 
jump processes
with rate $1$ and uniform initial distribution). 
To explain why the walks of \cite{BHPS03} can behave so differently
from those of the discrete web, we note that a single switch in the dynamical
random walk of~(\ref{other-dy}) affects 
only one increment of the walk while single
switches in the discrete web can change the path of the walker by a
``macroscopic'' amount. Indeed, the difference between the path 
$S_0^\tau$ before and after a single switch
is given by the difference between two independent simple 
random walks starting two spatial units apart.
This corresponds to the excursion of a (non-simple) random walk from zero whose mean duration
is infinite. It follows that a simple random walk is more sensitive to the extra 
noise induced by the dynamics on the discrete web than to the one induced by the dynamics considered in~\cite{BHPS03}. Rephrasing~\cite{SS05} in our context, since 
our dynamical random walk
``changes faster'' than the one in~\cite{BHPS03}, it 
has ``more chances'' to exhibit exceptional behavior. 
Mathematically,  ``changing fast'' corresponds to having
small correlations over short
time intervals and the
main ingredient for proving our exceptional times results
will be the correlation estimate~(\ref{speed-change})
of Proposition~\ref{sensitivityy-noise}.

By an obvious symmetry argument, there are also exceptional 
dynamical times $\tau$
for which $S^\tau(t) \leq j+K\sqrt{t}$ for all $t$. One may ask whether there are
exceptional $\tau$ for which $|S^\tau(t)| \leq j+K\sqrt{t}$ for all $t$.
Proposition~\ref{tme-K} below shows,
at least for small $K$, that there are no such exceptional times.
The case of large $K$ is unresolved.

{\bf Scaling Limits} 

In Section \ref{last-section2}, we discuss 
the continuum analogue of the dynamical random walk, 
the dynamical Brownian motion constructed in \cite{NRS08}.
We briefly recall there 
the main ideas of the construction 
along with some elementary properties of that object.
Then, we outline the main ideas that are needed to extend 
the results for exceptional times from the discrete level
to the continuum.

\section{Tameness}
\label{tameness}
In this section, we prove the following theorem.
\bteo
\label{th1} $\P\left(S^\tau_0 \ 
\txt{is recurrent for all $\tau\geq0$}\right)=1$.
\eteo
Recall the definition of the dynamical percolation model 
given in the introduction. 
For Bernoulli percolation on a homogeneous graph with
critical probability $p_c$, 
let $\theta(p)$ be the probability that
the origin belongs to an infinite cluster. In Section~3 of~\cite{HPS97}, 
it is proved that if
for some $C<\infty$
\begin{equation}
\theta(p)\leq C(p-p_c) \ \ \txt{for} \ \ p \geq p_c,
\end{equation}
then in the corresponding dynamical percolation model,
there is almost surely 
no dynamical time~$\tau$ at which percolation occurs. 
In our setting, an entirely parallel argument can be used to show
tameness of the dynamical discrete web with respect to
recurrence. 

Following~\cite{HPS97} we start by giving a very general tameness criterion. 
Let $\P_p$ be the probability measure for the static web when the 
probability for having a right arrow at a given site of $\Z_{even}^2$ 
is~$p$. Let $S_0$ be the 
simple random walk starting from the origin and 
let $A$ be 
a measurable
set of paths
such that $\P_{1/2}(S_0\in A)=0$.
In the following, we denote $\P_p(A)$ by $\theta_A(p)$.
Our first lemma is the analogue to
Lemma 3.1 in \cite{HPS97}.
\blem
\label{crit-tam}
Let $A$ be such that $\{S_0\in A\}$ is 
an increasing event w.r.t. the basic Bernoulli 
$\{\xi_z\}$ variables and such that $\P_{1/2}(S_0\in A)=0$.
Let $N_A$ be the cardinality of the set $\{\tau\in[0,1] \ :\  S^\tau_0\in A\}$.
Suppose there exists $c<\infty$ such that
\begin{equation}
\theta_A(p)\leq c(p-\frac{1}{2}) \ \ \txt{for all} \ \ p \geq \fr{1}{2},
\end{equation}
Then $\E(N_A)<\infty$.
\elem
\begin{proof}
Let $m>1$. We first estimate $\E(N_m)$ where 
$N_m$ is the number of $i\in\{1,2,...,m\}$
such that
there exists $\tau\in[\fr{i-1}{m},\fr{i}{m}]$ for which $S^\tau_0\in A$. 
For a given $i\leq m$, define
\begin{equation}
\bar \xi =\{\sup_{\tau\in[\fr{i-1}{m},\fr{i}{m}]} \xi_z(\tau)\}_{z\in\Z_{even}^2}.
\end{equation}
This naturally induces a new arrow configuration $\bar \aleph$ for which 
the probability
to find a right arrow at any given site is given by
\begin{equation}
\bar p=1-\frac{1}{2} \exp(-\frac{1}{2m})\leq \fr{1}{2} +\fr{1}{4m}.
\end{equation}
For such a configuration, the path $\bar S_0$ starting from the origin is 
a drifting random walk 
coupled with $S^\tau_0$ in such a way that
\begin{eqnarray}
\forall\tau\in [\fr{i-1}{m},\fr{i}{m}], \ \ S^\tau_0\leq \bar S_0,
\end{eqnarray}
which implies
\begin{eqnarray}
\P(\exists\tau\in[\fr{i-1}{m},\fr{i}{m}] \ \ \txt{with} \ \ \
S^\tau_0\in A) & \leq & \P(\bar S_0\in A) \\
& \leq &c(\bar p-\frac{1}{2})\leq\fr{c}{4m}.
\end{eqnarray}
Hence, $E(N_m)\leq m\frac{c}{4m}=c/4$ for all $m>1$.
Since $N_A= \liminf_{m\uaw\infty} N_m$, Fatou's lemma 
completes the proof.
\end{proof}

We now turn to the proof of Theorem \ref{th1}. For any $n\geq0$, 
let $A_n$ be the set of (piecewise linear) simple random walks
$\pi$ starting from the origin and such that for all $t\geq0$, $\pi(t)>-n$. 
It is well known that
\begin{equation}
\theta_{A_n}(p)=1-\left(\ 1-\fr{(2p-1)}{p}\ \right)^n, 
\ \ \ \ \txt{for}\ p\in[\fr{1}{2},1].
\end{equation}
Clearly, $A_n$ satisfies the hypotheses of Lemma \ref{crit-tam}, 
implying that $\E(N_{A_n})<\infty$.

In Lemma 3.2 of \cite{HPS97}, it is proved that
for any homogeneous graph with
critical probability~$p_c$, 
the number $N$ of times $\tau\in[0,1]$ such 
that in dynamical percolation the origin belongs
to an infinite cluster
is a.s. either $0$ or $\infty$.
By exactly the same reasoning, one can 
show that $N_{A_n}$ is either~$0$ or~$\infty$.
Since $\E(N_{A_n})<\infty$ 
for every $n$, 
we have that $N_{A_n}=0$ for every $n$ and this together
with the corresponding
result for transience to $-\infty$ completes the proof of Theorem~\ref{th1}.
 
 \bigskip

 \brm\label{grw}
 Another property of the static discrete web with
 respect to which the dynamical one is tame is the almost sure
 coalescence of all of its paths. Indeed, the difference between two 
 independent random walks 
 is 
 again a (non-simple) random walk and the proof of Theorem~\ref{th1}
 can easily be adapted to show that
 at every dynamical 
 time $\tau$ two walkers always meet and coalesce after some finite time~$t$.
 \erm

\section{Sensitivity to the Dynamics}
\label{sens}

In the following, $(C([0,1]),|.|_{\infty})$ denotes the space of 
continuous functions on $[0,1]$ equipped with the sup norm.
In order to prepare for our results about exceptional times, we need 
to prove that the arrow configuration in the DyDW
decorrelates fast enough to allow exceptional behavior for the dynamical 
random walk. This will be done by proving that on a large (diffusive) 
scale and for $\tau\neq\tau'$,  the paths $S^\tau_0$ and $S^{\tau'}_0$ 
evolve almost independently. More precisely, if for any (small) $\d>0$ 
and any $\pi\in C([0,1])$ we set  
$\tilde \pi(t)\equiv \pi(t/\d^2) \ \d$, we will prove that for a 
certain open set $O\in C([0,1])$, we have 
the following decorrelation inequality:
\beq
\label{noise-sensitivity2}
\P(\tilde S^\tau_0\in O \ , 
\ \tilde S^{\tau'}_0\in O )\leq \P(\tilde S\in O)^2 \ + 
\ K (\fr{\d}{|\tau-\tau'|})^a,
\eeq
where $S$ is a simple symmetric
random walk and $K$ does not depend on $\d,\tau$ and $\tau'$.
In other 
words, the inequality (\ref{noise-sensitivity2}) 
estimates the sensitivity of
the event $O$ to the dynamics. 

We now turn to our specific choice for $O$. 
Recall that we 
aim to prove that at some exceptional
$\tau$'s
the path $S^\tau_0$ is $K^{+}$-subdiffusive, which requires that the
walk starting from the origin is abnormally tilted
to the right.
Hence, it is natural to
study the noise sensitivity
of the event
\beq
\label{equtio-fr}
O=\{\forall t\in[0,1], \ \pi(t)>-1 \ \ \txt{and} \ \ \pi(1)>1\}
\eeq
which occurs for paths slightly tilted to the right. 
Studying noise sensitivity for this event is analogous to the 
corresponding question concerning left-right crossing of a square
in dynamical percolation as studied in \cite{SS05}.
The previous discussion
motivates the following proposition.

\bprop
\label{sensitivityy-noise}
For $O=\{\forall t\in[0,1], \ \pi(t)>-1 \ \ \txt{and} \ \ \pi(1)>1\}$,
there exist $K,a\in(0,\infty)$ (independent of $\d,\tau$ and $\tau'$) such that 
\beq
\label{speed-change}
\P(\tilde S^\tau_0\in O \ , \ 
\tilde S^{\tau'}_0\in O )\leq 
\P(\tilde S\in O)^2 \ + \ K (\fr{\d}{|\tau-\tau'|})^a,
\eeq
where $S$ is a simple symmetric random walk.
\eprop

In order to prove the proposition, we start by highlighting the fact that 
along the $t$-axis, the pair $(S_0^\tau,S_0^{\tau'})$ alternates between 
times  at which the two paths are equal (they ``stick together'') and 
times  at which they move independently. 
Recall that if $S^{\tau}_0$ and $S^{\tau'}_0$ coincide at time $t$ {\it and} if
the clock
at $(S^{\tau}_0(t),t)$ does not ring on $[\tau,\tau')$, the increments 
of $S^\tau_0$ and $S^{\tau'}_0$ at time $t$ are equal (i.e., $S^\tau_0$ 
and $S^{\tau'}_0$ stick together). Otherwise, the two increments are independent.
This suggests the following time decomposition of the pair 
$(S^{\tau}_0,S^{\tau'}_0)$.
Define inductively 
$\{T_k\}_{k\geq0}$ with $T_0=0$ and for any $k\geq0$,
\begin{eqnarray*}
T_{2k+1}& =& \inf\{n\in\N, n\geq T_{2k}: \  \txt{the clock at}\ (S^\tau_0(n),n)
\ \textrm{rings in} \  [\tau,\tau') \}, \\[0.1in]
T_{2k+2}& = & \inf\{n\in\N, n> T_{2k+1}: S^\tau_0(n)=S^{\tau'}_0(n) \}, \\[0.1in]
\D T_k & = & T_{2k+1}-T_{2k} \ \ \ \textrm{with $ \, \mathbb{P}( \D T_k \geq j) 
= e^{-|\tau-\tau'|j}$}.
\end{eqnarray*}
On the interval of integer time $[T_{2k},T_{2k+1}]$, the paths $S^\tau_0,S^{\tau'}_0$
coincide and at time $T_{2k+1}$ 
they move independently until meeting at time $T_{2k+2}$.
Hence, if we skip the intervals
$\{[T_{2k},T_{2k+1})\}_{k\geq0}$, $(S^\tau_0, S^{\tau'}_0)$ 
behave as two independent
random walks $(S^\tau_d,S^{\tau'}_d)$, while if we skip
$\{[T_{2k+1},T_{2k+2})\}_{k\geq0}$, the two walks coincide with a single
random walk $S_s$.
Furthermore, since $S_s$ is constructed from the arrow configuration at
different sites than the ones used to construct $(S_d^\tau,S_d^{\tau'})$, it is
independent of $(S_d^\tau,S_d^{\tau'})$. 

Now, skipping the intervals 
$\{[T_{2k},T_{2k+1})\}_{k\geq0}$ corresponds to
making the random time change
$t\rightarrow C(t)$  where $C$ is the right continuous inverse of 
\begin{eqnarray}
t+\sum_{k\leq l(t)} \D T_k \nn\\
\txt{with} \ \ l(t)=\#\{i\in\N, i\leq t: S^\tau_d(i)=S^{\tau'}_d(i)\}. \label{di-lo}
\end{eqnarray} 
Skipping $\{[T_{2k+1},T_{2k+2})\}_{k\geq0}$ corresponds to
making the time change $t\rightarrow t-C(t)$. This analysis yields 
the following lemma.

\blem
\label{pai-sbm}
There exist three independent simple symmetric
random walks $(S_s,S^{\tau}_d,S^{\tau'}_d)$ and an independent sequence
of independent non-negative integer valued random variables 
$\{\D T_k\}_{k\geq0}$ with $\P(\D T_k \geq j)=e^{-(\tau'-\tau)j}$
such that
\begin{eqnarray}
S^\tau_0(t)=S_d^{\tau}(C(t))+S_s(t-C(t)), \label{eq-stau} \\
S^{\tau'}_0(t)=S_d^{\tau'}(C(t))+S_s(t-C(t))\label{eq-stau'},
\end{eqnarray}
where $C$ is the right continuous inverse of (\ref{di-lo}).
\elem

In the following, the  pair $(S_0^\tau,S_0^{\tau'})$
will be referred to as a {\it sticky pair of random walks}. 
We note that the previous lemma has a continuous analogue called a 
sticky pair of Brownian motions---see Section \ref{last-section2} for more details.

Heuristically, in order to prove Proposition \ref{sensitivityy-noise}, 
we need to show that at large (diffusive) scales
Equations (\ref{eq-stau})-(\ref{eq-stau'}) become
\beqn
S_0^\tau(t)\approx S_d^\tau(t), \\
S_0^{\tau'}(t)\approx S_d^{\tau'}(t),
\eeqn
or equivalently that $C(t)\approx t$ (see Lemma \ref{lemm1} below).
The following three lemmas prepare
the justification of this informal approximation. Let $\d>0$. We recall
that for a path $S$, $\tilde S(\cdot)\equiv S(\cdot /\d^2)\ \d$. 
In the following, we set $\D\equiv \d/|\tau-\tau'|$ and for $O\subset C([0,1])$ 
and any $r\geq0$, we define 
$$O+r\equiv \{\pi\in C([0,1]) \ \txt{s.t.}
\ \exists \ \bar \pi\in O \  \txt{s.t.} \ |\pi-\bar \pi|_{\infty}\leq r \}.$$ 

\blem
\label{lo-u}
Let $S$ be a simple symmetric random walk.
For the $O$ defined in (\ref{equtio-fr}) and any $\a<\frac{1}{2}$,
\begin{equation}
\P(\tilde S\in [O+\D^\a]\setminus O) \leq \  c' \ \Delta^\alpha \label{lo-u1}
\end{equation}
where $c'\in(0,\infty)$ is independent of $\D$ and $\d$.
\elem
\begin{proof}
\beqn
\P(\tilde S\in [O+\D^\a]\setminus O) & \leq&  
\P(\inf_{t\in[0,1]} \tilde S(t)\in(-1-\D^\a,-1])+\P(\tilde S(1)\in(1-\D^\a,1]).\nn
\eeqn
We will prove that the second term on the right hand side of
of the inequality is of order $\D^\a$. The first term
can be handled similarly.

In~\cite{F73}, it is proved that a 
sequence of rescaled standard random walks 
$\{S(\cdot / \d^2) \d\}_{\d>0}$ 
and a Brownian motion $B$ can be constructed on the same probability
space in such way that for any $\a<\frac{1}{2}$ the quantity 
$\mathbb{P}(|B-S(\cdot \ / \d^2) \d|_{\infty}>\d^{\a})$ goes to $0$ faster than any 
power of~$\d$. On this probability space,
\begin{eqnarray}
\label{91}
&& \P(\tilde S(1)\in(1-\D^\a,1]) \nn \\
&& \leq   \P(B(1)\in [1-2\D^\a,1+\D^\a]) \  + 
\ \ {P}(|\tilde S-B|_\infty \geq \D^\a), 
\end{eqnarray}
Let $\a<\fr{1}{2}$. Because $\D^\a>\d^\a$ (since $|\tau-\tau'|\leq1$),  
the last term on the right-hand side
of (\ref{91}) is bounded by 
$O(\d)$, and consequently by $O(\D)$. By a density argument, the first term on
the right hand side of the inequality
is clearly bounded by $c \D^\a$ and this completes the proof of the lemma.
\end{proof}

\blem
\label{lemm1}
Define ${\bar C}(t)=C(t / \d^2) \d^2$ where $C$ is defined in (\ref{di-lo}) 
(note that the random clock $C$ 
and the paths are rescaled in a different manner).

For any $1>\beta>0$
\begin{equation}
\label{lemma3bound}
\mathbb{P}(\sup_{t\in[0,1]}\ (t-{\bar C}(t)) \geq {\Delta^\beta})\leq
\tilde c \Delta^{1-\beta},
\end{equation}
where $\tilde c\in(0,\infty)$ is independent of $\D$ and $\d$.
\elem

\begin{proof}
Recall from Lemma \ref{pai-sbm} that,
\begin{eqnarray}
\label{emmanuel}
C^{-1}(t)& =t+ L(t), \nn\\
\txt{where} \ \ \ L(t) &= \sum_{k\leq l(t)} \D T_{k},
\end{eqnarray}
$\{\D T_k\}$ are independent geometric random variables, $l$ is the 
discrete local time at the origin
of $S_d^\tau-S_d^{\tau'}$ and $C^{-1}$ is the right continuous inverse of $C$.
In the following, we set $\bar L(t)\equiv L(t / \d^2)\ \d^2$. We first prove that
\begin{equation}
\label{lemma3bound'}
\mathbb{P}(\bar L(1) \geq {\D^\b})\leq
\tilde c \D^{1-\b}.
\end{equation}

By the Markov inequality,
\begin{eqnarray}
\mathbb{P}({\bar L}(1)\geq{\D^{\b}})\leq 
\E\left({l(1/ \d^2)}{\d}\right) \ \ ({\E(\D T_1)}{\d}) \ \ \fr{1}{\D^\b}
\\[0.1in]
\textrm{with $\mathbb{E}(\D T_1)=\sum_{k=1}^\infty e^{-|\tau-\tau'| \ k}=
\frac{\exp(-|\tau-\tau'|)}{1-\exp(-|\tau-\tau'|)}$} \ .
\end{eqnarray}
Now
$$\E({l(1/\d^2)}) \, = \, \sum_{k \leq 1/ \d^2} 
\mathbb{P}(S_d^\tau (k) - S_d^{\tau'} (k) =0)$$
and it is a standard fact that the probability in the summation is $O(1/ \sqrt{k})$
as $k \to \infty$; thus
$\E({l(1/\d^2)} {\d})$ is uniformly bounded in $\d$ as $\d \to 0$. Furthermore,
since $E(\D T_1)=0(|\tau-\tau'|)^{-1}$, we have $\d E(\D T_1)=O(\D)$  
and thus 
(\ref{lemma3bound'})
follows.


Next, on the event
$\{{\bar L}(1) \leq {\Delta^{\beta}}\}$, 
(\ref{emmanuel}) implies that for any $t\in[0,1]$:
\begin{equation}
({\bar C})^{-1}(t)\leq t +\Delta^{\beta}.
\end{equation}
Since $\bar C(t)\leq t$ and ${\bar C}$ is an 
increasing function of $t$, it follows that
on $\{{\bar L}(1) \leq {\Delta^{\beta}}\}$,
for all $t\in[0,1]$, we have
\begin{equation}
t-\bar C(t) \leq \Delta^{\beta}.
\end{equation}
The lemma thus follows from (\ref{lemma3bound'}).
\end{proof}

\blem
\label{grr}
For any continuous function $f$, define 
$\omega_{f}(\e)=\sup_{s,t\in[0,1], |s-t|<\e} |f(t)-f(s)|$ to be
the modulus of continuity of $f$ on $[0,1]$.

Let $\a,\beta\in(0,\infty)$ be such that $\beta/2>\a$. 
For any $r\geq0$, there exists $c$ (independent of $\D$ and $\d$) such that
\beq
\P(\omega_{\tilde S}(\D^\b)\geq \fr{\D^\a}{2})\leq \ c \  \D^r \ .
\eeq
\elem
\begin{proof}
Let $m,n\geq0$
and define 
\beq \label{integformula}
\tilde M\equiv\int_0^1\int_0^1 \ \fr{|\tilde S(t)-\tilde S(s)|^n}
{|t-s|^m}\ dt ds.
\eeq
By the Garsia, Rodemich and Rumsey inequality~\cite{GRR70}, 
we have for $m>2$ and all $s,t\in[0,1]$
\beq
\label{gggri}
|\tilde S(t)-\tilde S(s)|\leq \fr{8m}{m-2} \ 
(4\tilde M)^{\fr{1}{n}} \ |t-s|^{\fr{m-2}{n}}.
\eeq
It is well known that $\E(|\tilde S(t)-\tilde S(s)|^n)
\leq c'|t-s|^{\fr{n}{2}}$, where $c'$ is uniform in $\d$.
Hence, (\ref{integformula}) implies that if $\fr{n}{2}-m>-1$,
then $\E(\tilde M)
\leq c<\infty$ so that for every $r\geq0$,
\beq
\label{m.s}
\P(\tilde M>\D^{-r})\leq c \D^r.
\eeq
On the other hand, on $\{\tilde M\leq \D^{-r}\}$, (\ref{gggri}) yields
\beqn
\omega_{\tilde S}(\D^\b)& \leq & \fr{8m}{m-2} 
\ (4 \D^{-r})^{\fr{1}{n}} \ |\D^\b|^{\fr{m-2}{n}} \\[0.1in]
& \leq &  c(n,m) \ \D^{\fr{1}{n}(\b(m-2)-r)}.
\eeqn
Since $\b/2>\a$, one can always take (for fixed $\a$, 
$\beta$, $r$) $n,m$ large enough such that
both
$\fr{n}{2}-m>-1$ and $\fr{1}{n}(\b(m-2)-r)>\a$.
For such a choice, and taking $\D$ small enough so that 
$c(n,m)\D^{\fr{1}{n}(\b(m-2)-r)}\leq \D^\a/2$, 
we obtain than on $\{\tilde M\leq \D^{-r}\}$
\beq
\omega_{\tilde S}(\D^\b)\ \leq \   \ \D^{\a}/2.
\eeq
Hence, for small enough values of $\D$, the 
claim of the lemma follows from (\ref{m.s}). The claim is
obviously satisfied for larger values of $\D$. 
\end{proof}

We are now ready to prove Proposition \ref{sensitivityy-noise}.
Recall the definition of $S_d^\tau$ and $S_d^{\tau'}$ in Lemma \ref{pai-sbm}. 
For any $\a>0$, we have
\beqn
\P(\tilde S^\tau_0\in O \ , \ \tilde S^{\tau'}_0\in O) & 
\leq \P (\tilde S^\tau_d\in O+\D^\a \ , \ \tilde S^{\tau'}_d \in O+\D^\a)
\nonumber \\[0.1in] 
& + 2\P(\tilde S^\tau_0\in O \ , \ \tilde S^\tau_d\in (O+\D^\a)^c) 
\label{equation1}. 
\eeqn
where $(O+\D^\a)^c$ is the complementary set of $O+\D^\a$.
Note that we used the equidistribution of
$(\sud,\su)$ and $(\sdd,\sd)$.
We start by
dealing with the first term on the right-hand side of the inequality. 
Since $\sud, \sdd$ are independent
and distributed like a rescaled simple symmetric 
random walk $\tilde S$, we have
\beqn
\P (\tilde S^\tau_d\in O+\D^\a \ , \ \tilde S^{\tau'}_d \in O+\D^\a)
& = \P(\tilde S^\tau_d\in O+\D^\a) \
\P(\tilde S^{\tau'}_d \in O+\D^\a) \nonumber \\
&\leq \P(\tilde S\in O)^2+ 2\P(\tilde S\in [O+\D^\a]\setminus O).\nonumber\\
\label{fab1}
\eeqn
The latter inequality and Lemma \ref{lo-u} above imply that 
\begin{eqnarray}
\label{au-suivant}
\P (\tilde S^\tau_d\in O+\D^\a \ , \ \tilde S^{\tau'}_d \in O+\D^\a)
\leq  \mathbb{P}(\tilde S\in O)^2 +
2 c' \D^{\a},
\end{eqnarray}
for any $\a<1/2$. By (\ref{equation1}) and (\ref{au-suivant}),  
Proposition~\ref{sensitivityy-noise} follows if there are $c'',a'\in(0,\infty)$
such that
\begin{eqnarray}
\P[\tilde S^\tau_0\in O \ , \ \tilde S^\tau_d\in (O+\D^\a)^c] \leq c'' \Delta^{a'}.
\end{eqnarray}
This inequality can be justified as follows.
Let $0<\beta<1$. By Lemma \ref{pai-sbm}
\begin{eqnarray}
\tilde S^\tau_0(t)& = \sud({\bar C}(t))+\Ss(t-{\bar C}(t)), \\[0.1in]
& = \sud(t)+ 
[\sud({\bar C}(t))-\sud(t)]+\Ss(t-{\bar C}(t)).
\end{eqnarray}
The last equality
implies that for any $0<\b<1$ with $\alpha<\beta / 2$,
\begin{eqnarray*}
\P[\tilde S^\tau_0\in O \ , \ \tilde S^\tau_d\in (O+\D^\a)^c]  
&\leq & \P(|\tilde S^\tau_0-\tilde S^\tau_d|_{\infty}\geq \D^\a) \\[0.1in]
& \leq &
\mathbb{P}(|\Ss (t-{\bar C}(t))|_{\infty} \geq \frac{\Delta^\alpha}{2})
+\mathbb{P}(|\sud(t)-\sud({ \bar C}(t))|_\infty 
\geq \frac{\Delta^\alpha}{2}) \\[0.1in]
& \leq &2 \mathbb{P}\left(\omega_{\tilde S}(\D^\b)
\geq \fr{\D^\a}{2}\right) + 2 \P(|t-\bar C(t)|_\infty \geq \D^\b) \\[0.1in]
&\leq  &2 c\D^r+2 \tilde c \D ^{1-\b}\, .
\end{eqnarray*}
where $r>0$ and the last inequality is given by  Lemmas 
\ref{lemm1} and \ref{grr} above.
So far we have only needed 
$\alpha\in(01/2)$   and thus we can indeed choose $\beta\in(0,1)$ 
and then $\a<\b/2$
so that
Proposition \ref{sensitivityy-noise} follows.

\section{Existence of Exceptional Times}
\label{existence}
In this section we prove the following result.
\bprop
\label{violation}
For $K$ large enough 
\begin{equation}
\P(\exists \tau\in[0,1], \ \ \txt{s.t. $S_0^\tau$ is $K^{+}$-subdiffusive})>0.
\end{equation}
(For a definition of $K^{+}$-subdiffusivity, see Definition \ref{K-subd}.)
\eprop

\bigskip

\begin{figure}
\centering
\includegraphics[scale=0.1333]{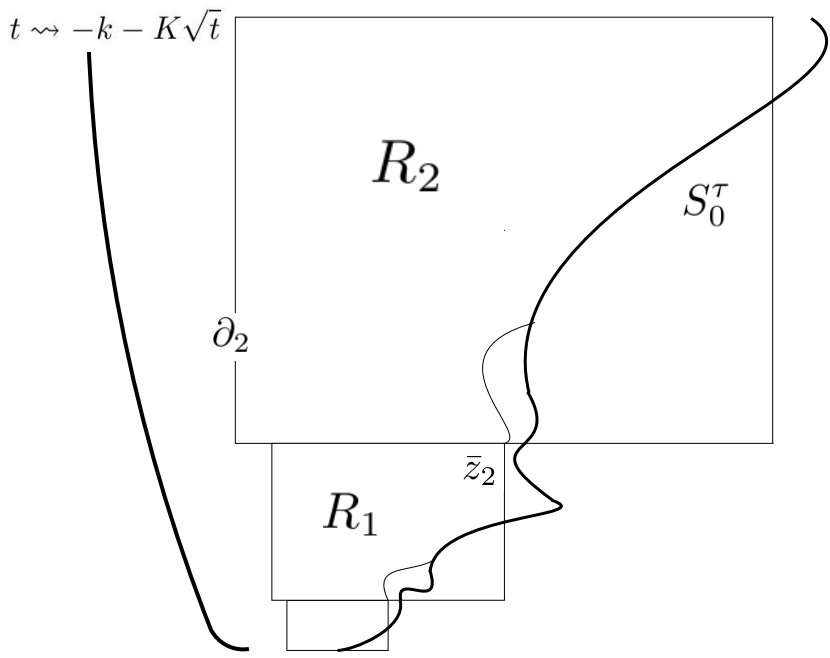}
\caption{
Construction of the first three boxes $(R_0,R_1,R_2)$
with $t$ the vertical and $x$ the horizontal coordinate. The thin curves
represent segments of the paths starting from 
$\bar z_i$, for $i=1,2,3,$ for which
the events $A_0^\tau$, $A_1^\tau$ and $A_2^\tau$ occur.}
\label{explain}
\end{figure}

Let $\g>2$ and $d_k=2(\lfloor \frac{\gamma^k}{2} \rfloor +1)$, 
where $\lfloor x \rfloor$ is the integer part
of $x$.
We construct inductively a sequence of ``diffusive'' boxes $R_k$ in the following manner (s
ee Figure \ref{explain}).
\begin{itemize}
\item $R_0$ is the rectangle with vertices $(-d_0,0)$,
$(+{d_0},0)$, $(-{d_0},d_0^2)$ and
$(+{d_0},d_0^2)$.
\item Let $\bar z_{k}=(x_k,t_k)$ be the middle
point of the lower edge of $R_k$ (e.g., $\bar z_{0}=(0,0)$). 
$R_{k+1}$
is the rectangle of height $d_{k+1}^2$ and width $2 d_{k+1}$ such that $\bar z_{k+1}$
coincides with the the upper right vertex of $R_k$ (see Figure~\ref{explain}).
\end{itemize}
Note that for our particular choice of $d_k$, $\bar z_k$ always belongs to 
$\Z_{even}^2$ for $k\geq0$ and a simple computation leads to the following lemma. 

\blem
\label{shape;curev}
Let 
${\partial}^{\g}={\partial}^{\g} (t)$ 
denote the right-continuous function obtained by joining
together the left boundaries ${\partial}_k$ of ${R}_k$. 
For any $K>0$, let $\g(K)$ be the solution in $(2,\infty)$ of 
$K=(\g-2)\sqrt{\frac{\gamma+1}{\gamma-1}}$. Then,
\begin{equation}
\label{g-ine}
\forall t\geq 0,\ \  \partial^{\g(K)} (t) \geq -3-K\sqrt{t}. 
\end{equation}
\elem
\begin{proof}
On $[t_n,t_{n+1})$, we have
${\partial}^\g(t) = {\partial}^\g({t}_n) =
{x}_n - {d}_n$ $=
(d_0+d_1+\dots+d_{n-1}- d_{n})$. 
If $\g$ is such that
\begin{equation}
\label{boundary}
{\partial}^\g({t}_n) \ \geq \ - (3+K \sqrt{{ t}_n}) \
\textrm{for} \  n = 0, 1, 2, \dots \ ,
\end{equation}
then we will have ${\partial}^\g(t) \geq -( 3+K \sqrt{t})$
for all $t \geq 0$ as desired. 

The inequality (\ref{boundary}) can be rewritten as 
\begin{equation}
\label{gggg}
d_n\leq 3 + d_0+\cdots+
d_{n-1}+ K[{d_0}^2+\cdots+d_{n-1}^2]^{1/2}. 
\end{equation}
Using the bound $d_{n}\leq 2+\gamma^n$ on the 
left-hand side of (\ref{gggg}) and the bounds 
$d_j\geq\gamma^j$ on the right-hand 
side, it follows that in order to verify~(\ref{gggg})
it suffices to have, for $n=0,1,2,\dots$,
\begin{equation}
\label{anotherineq}
\gamma^n\leq 1 +
\frac{\g^n-1}{\g-1}+K\sqrt{\frac{\gamma^{2n}-1}{\gamma^2-1}}\ .
\end{equation}
Using the elementary bound
$\sqrt{\gamma^{2n} - 1} \geq \gamma^n(1 - \gamma^{-2n})$ (for $\gamma \geq 1$),
we see that in order to verify~(\ref{anotherineq}),
it suffices to have, for $n=0,1,2,\dots$,
\begin{equation}
\gamma^n(\frac{\gamma-2}{\gamma-1} - \frac{K}{\sqrt{\gamma^2 -1}})
\leq 1 -
\frac{1}{\gamma-1}- \frac{K}{\sqrt{\gamma^2-1}}\gamma^{-n} \, .
\end{equation}
Choosing $\g$ such that $K=({\gamma-2}) \sqrt{\frac{\gamma+1}{\gamma-1}}$ 
yields (\ref{g-ine}). The lemma follows from the fact that 
$\g\raw({\gamma-2}) \sqrt{\frac{\gamma+1}{\gamma-1}}$ 
is a continuous increasing function mapping $(2,\infty)$ onto $(0,\infty)$.
\end{proof}
By Lemma \ref{shape;curev}, $S_0^\tau$ is $K^{+}$-subdiffusive if 
$S^{\tau}_0(t)\geq \partial^{\g(K)}(t)$.
Let $S^\tau_{\bar z_k}$ be the path in $W(\tau)$ starting 
from $\bar z_k=( x_k, t_k)$ and define the event
\begin{equation}
\label{def-akt}
A_k^\tau = 
A_k^\tau(K)=\{\forall t\in [t_{k},t_{k+1}] \ \ S^\tau_{\bar z_k}(t)> \partial_k(t) \ 
\  ,  \  S^\tau_{\bar z_k}(t_{k+1}) > x_{k+1}\}.
\end{equation}
(Here $\partial_k$ depends implicitly on $\g(K)$.)
Since paths in $W(\tau)$ do
not cross, if $\cap_{k\leq n}A_k^\tau$ occurs,  
$S_0^\tau$ is forced to remain to the right of $\partial_k$ 
on $[t_k,t_{k+1}]$ for every $k\leq n$ (see Figure \ref{explain}). This implies that 
if we have
\begin{equation}
\label{0ba-oba}
\P(\exists \tau\in[0,1], \ \ \cap_{k\geq 0} A_k^\tau(K) \ \ \txt{occurs})>0,
\end{equation}
then
\begin{equation}
\P(\exists \tau\in[0,1], \ \ S_0^\tau  \ \ \txt{is $K^{+}$-subdiffusive})>0.
\end{equation}
In the rest of the section we proceed to verify 
(\ref{0ba-oba}).

In the following, $K$ is temporarily fixed and to ease the notation we 
write $A_k^\tau$ for $A_k^\tau(K)$ and $\g$ for $\g(K)$. In order to 
verify (\ref{0ba-oba}), we start by proving the following lemma using 
Proposition~\ref{sensitivityy-noise}.
\blem\label{jrn} There exists $c\in(0,\infty)$ such that for $\tau,\tau'\in[0,1]$
\begin{eqnarray}
\label{energy}
 \forall n\geq 0, \ \ \prod_{k=0}^n
 \frac{\mathbb{P}(A_k^\tau \cap A_k^{\tau'})}{\mathbb{P}(A_k)^2} 
 \leq {c} \frac{1}{|\tau-\tau'|^b},
 \end{eqnarray}
 where $A_k\equiv A_k^0$ and $b=\log(\sup_k [\P(A_k)^{-1}])/\log \g>0$.
 \elem
 \begin{proof}
 Let $(S^\tau,S^{\tau'})$ be the paths starting at $(0,0)$, defined as the
 translated version of the
 pair $(S^\tau_{\bar z_k},S^{\tau'}_{\bar z_k})\in(W(\tau),W(\tau'))$ 
 starting at $\bar z_k$. By translation invariance, $(S^{\tau},S^{\tau'})$
 is a sticky pair of random walks starting at $(0,0)$
 whose  distribution is 
described in Lemma \ref{pai-sbm} and by definition
 \beq
 A_k^\tau=\{S^\tau(d_k^2)>d_k\ \ ,\  \ \inf_{[0,d_k^2]}S^\tau(t)>-d_k\}.
 \eeq
 By Proposition \ref{sensitivityy-noise} for $\d=d_k^{-1}$,
 there exists $c,a\in(0,\infty)$ such that
 \beq
 \P(A_k^\tau \ , A_k^{\tau'}) \leq \P(A_k)^2 + c \, (\fr{1}{\g^k|\tau-\tau'|})^a.
 \eeq 

Defining $N_0=[\frac{-\log(|\tau-\tau'|)}{\log{\ \gamma}}]+1$
so that $(\g^{N_0} |\tau-\tau'|)\geq 1$, we have for $n>N_0$
\begin{eqnarray}
\label{final1}
\prod_{k=N_0+1}^{n} \left(
\frac{\mathbb{P}(A_k^\tau \cap A_k^{\tau'})}{\mathbb{P}(A_k)^2} \right) & \leq &
\prod_{k=N_0+1}^\infty (1+ \frac{c / \mathbb{P}(A_k)^2 }
{ |\tau-\tau'|^a\gamma^{a N_0} \ \ \gamma^{a(k-N_0)}}), \nonumber \\
& \leq &\prod_{k=1}^{\infty} (1+ \frac{c}{\inf _n\mathbb{P}(A_n)^2}
\ \ \frac{1}{\g^{ak}}). \label{final3}
\end{eqnarray}
The right-hand side of
(\ref{final1}) is independent of
$|\tau-\tau'|$ and is finite. Indeed, we have $0<\inf_n P(A_n)$  since the
boxes $R_k$ have diffusively scaled  sizes and therefore 
$\mathbb{P}(A_k)\rightarrow\mathbb{P}(A)$ as $k \to \infty$, where $A$ is the
event that a Brownian motion
${B}(t)$ starting at $0$ at
time $0$ has ${B}(1)>1$ and $\inf_{t\in[0,1]} {B}(t)>-1$.

On the other hand, for $n\leq N_0$
\begin{eqnarray}
\label{final2}
\prod_{k=0}^{n} \frac{\mathbb{P}(A_k^\tau \cap A_k^{\tau'})}
{\mathbb{P}(A_k)^2} & \leq (\sup_k \frac{1}{\mathbb{P}(A_k)})^{N_0+1},\nonumber \\
& \leq  c''\exp(\frac{\log[\sup_k
(\P(A_k)^{-1})]}{\log{\gamma}}\log(\frac{1}{|\tau-\tau'|})), \nonumber \\[0.1in]
& =  {c''}/{|\tau-\tau'|^b},  
\end{eqnarray}
where $c''=\sup_k (\P(A_k)^{-1})$ and 
$b=\log[\sup_k (\P(A_k)^{-1})] / \log{\g}$ are in $(0,\infty)$.
This and (\ref{final3}) imply (\ref{energy}).
\end{proof}

Following ~\cite{SS05}, the
Cauchy-Schwarz inequality and the previous lemma imply that
for every $n\geq0$, we have
\begin{eqnarray}
\mathbb{P}(\int_0^1 \prod_{k=0}^n 1_{A_k^\tau} d\tau>0) & \geq &
\frac{\left(\mathbb{E}\left[\int_0^1\prod_{k=0}^n 1_{A_k^\tau} 
\ d\tau\right]\right)^2}
{\mathbb{E}\left[\left(\int_0^1 \prod_{k=0}^n 1_{A_k^\tau} 
\ d\tau \right)^2\right]}, \\[0.1in]
& = &\left(\left[\int_0^1 \int_0^1 \prod_{k=0}^n
\frac{\mathbb{P}(A_k^\tau \bigcap A_k^{\tau'})}
{\mathbb{P}(A_k)^2} \ d\tau \ d\tau' \right]
\right)^{-1}\\[0.1in]
&\geq & {c^{-1}} \left(\left[\int_0^1 \int_0^1 
\frac{1}{|\tau-\tau'|^b} \ d\tau \ d\tau' \right]
\right)^{-1}\label{integrand}
\end{eqnarray}
where the equality is a consequence of the stationarity
of $\tau\rightarrow W(\tau)$ and the independence between the arrow configurations 
in different boxes $R_k$.
Recall that $\g$ has an implicit dependence on $K$ and that
$\g$ increases from $0$ to $\infty$ as $K$ increases on 
$(0,\infty)$ (see Lemma \ref{shape;curev}).
Hence, for K large enough 
such that $\g=\g(K)>\sup_k {\mathbb{P}(A_k)}^{-1}$, we have 
$$b=\log(\sup_k [\P(A_k)^{-1}])/\log \g(K)<1$$ and
$(\tau,\tau')\rightarrow |\tau-\tau'|^{-b} \in 
L^1([0,1]\times[0,1], d\tau \ d\tau')$. (\ref{integrand}) then
implies that
\begin{equation}
\label{esd}
\inf_n \mathbb{P}(\int_0^1 \prod_{k=0}^n 1_{A_k^\tau} \ d\tau > \ 0 ) \geq p>0.
\end{equation}
Let $E_n$ be the set of times $\tau$ in $[0,1]$ such that
$\bigcap_{k=0}^n A_k^\tau$ occurs. (\ref{esd}) implies that
$\mathbb{P}(\bigcap_{n=0}^{\infty}\{E_n\neq \emptyset\}) \geq p >0$.
Since $\{E_n\}$
is obviously decreasing in $n$, if the $E_n$ were closed
subsets of $[0,1]$ it would follow that $\mathbb{P}((\bigcap_{n=0}^{\infty} E_n)
\  \neq \emptyset) \geq p > 0$.

Unfortunately, the set of times at which one arrow is
(or any finitely many are) oriented to the right (resp., to the left)
is not 
in general 
a closed subset of $[0,1]$ since we have a right continuous process,
and thus $E_n$ is not in general a closed set. 
This extra technicality is
handled like in Lemma 5.1 in~\cite{SS05},
as follows. On the one hand, there are only countably many switching times
for all $\xi_{z}^\tau$'s (recall that $\xi_z^\tau$ represents the arrow
direction at location $z$).
On the other hand, at any switching time 
$\tau$, $\, \bigcap_{n\geq0} A_n^\tau$
does not occur  
by independence of the $\xi_{z}^\tau$'s. Since there are 
countably many switching times,
this implies that almost surely, the closures $\bar E_n$
of $E_n$ satisfy
\begin{equation}
\label{clos-ex}
\cap_{n=1}^\infty {\bar E}_n = \cap_{n=1}^\infty E_n.
\end{equation}
This completes the verification of (\ref{0ba-oba}) and thus the proof of
Proposition~\ref{violation}.

\section{Hausdorff Dimension Of Exceptional Times}
\label{hausdorff}
In this section, we derive some lower and upper bounds for the Hausdorff dimension of
the set of
exceptional dynamical times $\tau \in [0,\infty)$ at which $S^\tau_0$ 
becomes subdiffusive. 

\bdf
\label{def1}
We say that $\tau$ is a $K$-exceptional time if the path
$S^\tau_0$ in $W(\tau)$ 
does not cross the moving boundary $t\rightsquigarrow-K\sqrt{t}$.
$\mathcal{T}(K)$ is then defined as the set of all $K$-exceptional times
$\tau\in [0,\infty)$.
\edf

Clearly, the set consisting of all the $K$-exceptional times in $[0,\infty)$
is a non-decreasing function of $K$. The next propostion asserts that 
for fixed $K$, the Hausdorff
dimension $dim_H$ of the set of exceptional times is unchanged if $-K\sqrt{t}$
is replaced
by $-j-K\sqrt{t}$ for any $j \geq 0$.
We note that as in dynamical percolation (see Sec.~6 of~\cite{HPS97}),
$dim_H(\mathcal{T}(K))$ is a.s. a constant by the ergodicity in $\tau$ of
the dynamical discrete web.

\bprop
\label{constant}
The Hausdorff dimension $dim_H$ of the set $T_j = T_j(K)$
of exceptional times $\tau \in(0,\infty)$ such
that $S^\tau_0$ does not cross the moving boundary $t\rightsquigarrow-j-K\sqrt{t}$
does not depend on $j \geq 0$ (for fixed $K$).
\eprop

\begin{proof}
By monotonicity in $j$, it is enough to prove that
$dim_H(T_{j})\leq dim_H(T_{0})$ for $j$ any 
positive integer.

First,
$T_0 \supset T_{j}'\bigcap
\{\tau\in [0,1]:\xi_{(m,m)}^\tau =+1 \ \ \textrm{for} \ \ m<j\}$ 
where $T_{j}'$ is the set of $\tau\in [0,\infty)$ 
such that $S^\tau_{(j,j)}(n)\geq-K\sqrt{n}$ 
for $n\geq j$. 
Furthermore, $T_{j}' \supset \bar T_{j}$, 
where $\bar T_{j}$ is the 
set of $\tau\in [0,\infty)$ such that 
$S^\tau_{(j,j)}(n)-j\geq-j -K \sqrt{n-j}$ for $n\geq j$. 
Note that $\bar T_{j}$ is just the translation (from $(0,0)$ to $(j,j)$) of 
$T_{j}$. Hence, $\dim_{H}(T_0)$ is 
at least the dimension of
$\{\tau\in[0,1]: \forall 
\  m<j, \ \xi_{(m,m)}^\tau=+1\} \bigcap \bar T_j$.

By ergodicity in $\tau$, the a.s. constant $dim_H(\bar T_{j})$
is the essential supremum of the random variable 
$dim_H(\bar T_j\bigcap[0,1])$. On the other hand, since 
$\bar T_j\bigcap[0,1]$ and $\{\tau\in[0,1]: \forall \  m<j, \ \xi_{(m,m)}^\tau=+1\}$ 
are independent and the probability to have 
$\{\forall \ \tau\in[0,1], \ \forall \  m<j, \ \xi_{(m,m)}^\tau=+1\}$ is strictly
positive, it follows that 
$dim_H(\{\tau \in [0,1]: \forall \  m<j, \ \xi_{(m,m)}^\tau=+1\} \bigcap \bar T_j)$ 
has the same essential sup as $dim_H(\bar T_j \bigcap [0,1])$. Hence
$dim_H(T_{j}) = dim_H(\bar T_j) \leq dim_H(T_{0})$
and the conclusion follows.

\end{proof}

\subsection{Lower Bound}
\label{lowerbound}
Set ${\gamma}_0 \equiv \sup_{k,K} 1/{\mathbb{P}(A_k(K))}$, where 
$A_{k}(K)$ is defined by (\ref{def-akt}) with $\tau=0$. (Note that $\g_0>2$.)
We recall that ${\g}(K)$ is the solution in $(2,\infty)$
of $K=K(\g)=({\g -2})\sqrt{\frac{\g+1}{\g-1}}$
for $K>0$.
In this section, we prove the following proposition using Lemma \ref{jrn} and 
then
arguments identical to ones in \cite{SS05}.
\bprop
\label{lower}
\begin{equation}
\label{lowerdim}
dim_H({\mathcal T}(K)) \ \geq \ 1\, - \,
\frac{\log{\gamma_0}}{\log{\gamma(K)}}
\ \textrm{for} \ K\, > \, K(\g_0) \, .
\end{equation}
Thus, $\lim_{K\uaw\infty} dim_H({\mathcal T}(K))=1$.
\eprop
Let $K>K(\g_{0})$. Note that since $K\raw\g(K)$ is increasing, $\g(K)>\g_0$.
In the following and as in Section \ref{existence}, we drop the dependence
on $K$ in the notation. 
Consider the random measure $\sigma_n$, such that for any Borel set $E$ in $[0,1]$
$$\sigma_n(E)=\int_E \prod_{k=0}^n \fr{1_{A_k^\tau}}{\mathbb{P}(A_k)} \ d\tau.$$
We note that $\sigma_n$ is supported by $\bar E_n$, the closure of $E_n$ with
\begin{equation}
E_n=\{\tau\in[0,1]: \cap_{k\leq n} A_k^\tau \ \ \txt{occurs} \}.
\end{equation} 
For any positive measure $\s$, define the $\alpha$-energy of $\sigma$
as
\begin{equation}
\mathcal{E}_{\alpha}(\sigma)=\int_{0}^1 \int_0^1 \frac{1}{|\tau-\tau'|^{\alpha}} \
d \sigma(\tau) \ d \sigma(\tau') \, .
\end{equation}
Following \cite{SS05}, we will need the following extension of Frostman's Lemma.
\blem{\bf\cite{SS05}}
\label{lower-ss05}
Let $D_1\supset D_2\supset...$ be a decreasing sequence of compact subsets of $[0,1]$,
and let $\mu_1,\mu_2,...$ be a sequence of positive measures with $\mu_n$
supported on $D_n$. Suppose that there exists a constant $C\in(0,\infty)$ 
and $\a\in(0,1)$ such that for infinitely many
values of $n$,
\begin{equation}
\mu_n([0,1])\geq 1/C, \ \ \ \mathcal{E}_{\a}(\mu_n)\leq C.
\end{equation}
Then the Hausdorff dimension of $\cap D_n$ is at least $\a$.
\elem

Using the ergodicity of the dynamical web in the variable
$\tau$, we will prove Proposition~\ref{lower}
by showing that for $\a<1-\fr{\log(\g_0)}{\log(\g(K))}$,
$\{\s_n\}$ satisfies the hypotheses of this lemma with strictly positive 
probability. By Lemma \ref{jrn}, we have for all $n$ that
\begin{eqnarray}
\label{had}
\E[\sigma_n([0,1])^2]
& =   \int_0^1 \int_0^1 \prod_{k=0}^n
\frac{\mathbb{P}(A_k^\tau\cap A_k^{\tau'})}{\mathbb{P}(A_k)^2} 
d \tau d\tau ' \nonumber \\[0.1in]
& \leq  c \left(\left[\int_0^1 \int_0^1
\frac{1}{|\tau-\tau'|^{b}} \ d \tau \ d \tau' \right] \right), \nonumber 
\end{eqnarray}
where 
\beqn
\label{avant-ca}
b=\log{[\sup_{k} (\P(A_k)^{-1})]}/\log{\gamma}\leq  \fr{\log(\g_0)}{\log{(\g)}}
<1.
\eeqn
By the Cauchy-Schwarz inequality
\beqn
\E\left[\sigma_n([0,1])^2\right]^{\fr{1}{2}} \ 
\P\left[\sigma_n({[0,1]})\ > \ \frac{1}{2}\right]^{\fr{1}{2}}
& \geq & 
\E\left[\sigma_n([0,1])\cdot 1_{\s_n([0,1])>1/2}\right] \nn \\
& \geq &\E\left[\sigma_n([0,1])\right]-\fr{1}{2} \ =\ \fr{1}{2} \nn,
\eeqn
which implies that $\P[\sigma_n([0,1]>\fr{1}{2})]>c_1$ for 
some $c_1>0$ not depending on $n$.

\bigskip

By Fubini's Theorem and Lemma \ref{jrn},
\begin{eqnarray}
\label{en2}
\mathbb{E}(\mathcal{E}_{\alpha}(\sigma_n))& =
\int_0^1 \int_0^1 |\tau-\tau'|^{-\a} \ \prod_{k=0}^n
\frac{\mathbb{P}(A_k^\tau\cap A_k^{\tau'})}{\mathbb{P}(A_k)^2}
\ \ d\tau \ d\tau'. \nonumber \\[0.1in]
& \leq \ c \ \int_0^1 \int_0^1
\frac{1}{|\tau-\tau'|^{b+\a}} \ d\tau \ d\tau' .
\end{eqnarray}
Taking $\a$ such that
\begin{equation}
\label{rel-a}
\a<1-\fr{\log(\g_0)}{\log{(\g)}},
\end{equation}
we have from (\ref{avant-ca}) that
$b+\a<1$ and therefore
\begin{eqnarray}
\label{Haus}
\sup_{n\geq0} \mathbb{E}(\mathcal{E}_{\alpha}(\sigma_n)) 
\leq c_2<\infty. 
\end{eqnarray}
By Markov's inequality, for all $n$ and all $T$,
\begin{equation}
\P(\mathcal{E}_\a(\s_n)\geq c_2 T)\leq 1/T.
\end{equation}
Choose $T$ such that $1/T<c_1/2$. Letting 
\begin{equation}
U_n^\a=\{\s_n([0,1])>\fr{1}{2}\}\cap \{\mathcal{E}_\a(\s_n)\leq c_2 T\},
\end{equation}
by the choice of $T$, we have that
\begin{equation}
\P(U_n^\a)\geq c_1/2.
\end{equation}
By Fatou's lemma,
\begin{equation}
\P(\limsup_{n\uaw\infty} U_n^\a)\geq c_1/2.
\end{equation}

By Lemma \ref{lower-ss05},
it follows that for $\a$ satisfying (\ref{rel-a}), $\cap_{n\geq0} \bar E_n$ 
has Hausdorff dimension at least~$\a$ with positive probability.
Since $\cap_{n\geq0} \bar E_n=\cap_{n\geq0} E_n$ (see (\ref{clos-ex})), 
the same statement holds
for  $\cap_{n\geq0} E_n$ and we are done.

\subsection{Upper Bound}
\label{upperbound}
We will prove the following proposition.
\bprop
\label{upper}
$\textrm{dim}_H(\mathcal{T}(K)) \leq 1-p({K})$
where $p({K}) \in (0,1)$ is the solution of the equation
\begin{equation}
\label{sato}
f(p,{K}) \equiv
\frac{\sin(\pi p / 2) \Gamma(1+ p / 2)}{\pi} \sum_{n=1}^{\infty}
\frac{(\sqrt{2} K)^n}{n !} \ \Gamma((n-p)/2)=1 \, .
\end{equation}
Furthermore, $K\rightsquigarrow p(K)$ is a continuous decreasing function
on $(0,\infty)$ with
\begin{equation}
\label{hygeo}
\lim_{K\uparrow\infty}p(K)=0 \ \ \ \txt{and
more significantly} \ \
\lim_{K\downarrow0}p(K)=1.
\end{equation}
\eprop

To prove Proposition~\ref{upper}
we need the following lemma proved in the Appendix.

\blem
\label{principal22}
Let $0<l<1$. Let $S_\e$  be the simple asymmetric random walk
with 
\beq\label{value_of_drift}
\mathbb{P}(S_{\epsilon}(n+1)-S_{\epsilon}(n)=+1)=
\frac{1}{2}+\frac{1}{2}(1-e^{-\epsilon}).\eeq
Then there exists $c(l)$ such that
\begin{equation}
\P(\forall n, \ S_\e(n)\geq-1-K\sqrt{n})\leq c(l) \e^{p(K/l)}
\end{equation}
where $p({K})$ is the real solution in $(0,1)$ of~(\ref{sato}) (which satisfies
(\ref{hygeo})).
\elem

Let us partition $[0,1]$ into intervals of equal length
$2\epsilon$, and select the intervals containing a $K$-exceptional time. The
union of those is a cover of $\mathcal{T}(K)$ and we now estimate
the number $n(\epsilon)$ of intervals in the cover.

Let $U_{\epsilon}$ be the event that there is a time $\tau$ in $[0,2\epsilon]$
such that  $\tau\in \mathcal{T}(K)$. From the full dynamical arrow configuration
for all $ \tau \in [0,2\epsilon]$, we construct a static arrow configuration
as follows. We declare 
the static arrow at $(i,j)$ to be 
right oriented if and only if the dynamical arrow is right oriented
(i.e., $\xi_{i,j}^\tau = +1$) at some 
$\tau\in[0,2\epsilon]$ (a similar construction was used in Section~\ref{tameness}).
In this configuration, the path $S_{\epsilon}$ starting
from the origin and following the arrows is a slightly right-drifting
random walk with 
$\mathbb{P}(S_{\epsilon}(n+1)-S_{\epsilon}(n)=+1)=
\frac{1}{2}+\frac{1}{2}(1-e^{-\epsilon})$. Clearly,
\begin{equation}
\mathbb{P}(U_{\epsilon})\leq
\mathbb{P}(\forall n, \, S_{\epsilon}(n)\geq -1-K \sqrt{n}) \, .
\end{equation}
Lemma~\ref{principal22} 
implies that for any $l<1$
\begin{equation}
\mathbb{P}(U_{\epsilon})
\ \leq c(l) \epsilon^{p(\frac{K}{l})}.\label{final-touch}
\end{equation}
Hence
\begin{equation}
\mathbb{E}(n(\epsilon)) = O ( \epsilon^{p(\frac{K}{l})-1}) 
\end{equation}
so that
\begin{equation}
\limsup_{\epsilon \to 0} \mathbb{E}(\frac{n(\epsilon)}
{\epsilon^{p(K/l)-1}}) \ < \ \infty \, .
\end{equation}
By Fatou's Lemma,
$\liminf_{\epsilon \daw 0} n(\epsilon) \,
\epsilon^{1-p(K/l)}$
is almost surely bounded, which implies that 
$dim_{{H}}\mathcal{T}(K)$ (which is equal to $dim_{{H}}(T_1(K))$ by 
Proposition~\ref{constant}) is bounded above by $1-p(\frac{K}{l})$ 
for any $l<1$. Since $p(K)$ is continuous in $K$,
Proposition~\ref{upper} follows. 

\brm
We conjecture that $1-p(K)$ is the exact Hausdorff dimension of $\mathcal{T}(K)$.
\erm

Finally, Lemma \ref{principal22} also yields the following tameness result.
\bprop
\label{tme-K}
Let $K_1,K_2>0$ be small enough so that
$p(K_1)+p(K_2)>1$, where $p(K)$ is defined in Proposition \ref{upper}. 
For any $j\geq0$,
\begin{equation}
\P(\exists \tau\in[0,1] \ \txt{s.t. $\forall t\geq0$,} \ -j-K_1\sqrt{t} 
\leq \ S^\tau_0(t) \ \leq +j+K_2\sqrt{t})=0.
\end{equation}
\eprop
\begin{proof}

Define $U_\e^{+}$ (resp., $U_\e^{-}$) to be the event that for some 
$\tau\in[0,2\epsilon]$
and all~$t\geq0$,
$S^\tau(t) \leq +k+K_2\sqrt{t}$ (resp., $S^\tau(t) \geq -k-K_1\sqrt{t}$).
$U_\epsilon^+$ (resp., $U_\epsilon^-$) is a decreasing
(resp., increasing) event with respect to the basic $\xi_{(i,j)}^\tau$
processes. Hence, using the
FKG inequality, we have 
$$\mathbb{P} (U_\epsilon^+ \cap U_\epsilon^-) \leq \mathbb{P} (U_\epsilon^+)
\, \cdot \, \mathbb{P} (U_\epsilon^-).$$
Reasoning as in Proposition \ref{upper}, for any $l<1$, we have
\begin{eqnarray}
\P(U_\e^{-})&\leq &\P(\forall t\geq0, \ \ S_{\e}(t)\geq -j-K_1\sqrt{t})\\[0.1in]
& \leq & c_1 \epsilon^{p(\frac{K_1}{l})},
\end{eqnarray}
where $S_\e$ is defined as in the proof of Proposition \ref{upper}. The
second inequality is given by (\ref{final-touch}) immediately for 
$j\leq1$ and with a little bit of extra effort for all $j$. Symmetrically,
\begin{equation}
\P(U_\e^{+})\leq c_2 \epsilon^{p(\frac{K_2}{l})},
\end{equation}
which implies that
\begin{equation}
\mathbb{P} (U_\epsilon^+ \cap U_\epsilon^-) \leq c_1 c_2 \ \e^{p(K_1/l)+p(K_2/l)}.
\end{equation}
Take $l$ close enough to $1$ so that $p(K_1/l)+p(K_2/l)>1$ and
define $N$ as the cardinality of $\{\exists\tau\in[0,1] \ \txt{s.t. 
$\forall t\geq 0$,} \ -j-K_1\sqrt{t} \leq \ S^\tau_0(t) \ 
\leq +j+K_2\sqrt{t}\}$. Reasoning as in Lemma \ref{crit-tam}, we have
\begin{equation}
\E(N) = \lim_{\e\daw0} \fr{1}{2\e} \P(U_\e^+ \cap U_\e^-)=0,
\end{equation}
which completes the proof of the proposition.
\end{proof}

\section{Scaling Limit}
\label{last-section2}
In this section, we discuss 
the existence of a dynamical Brownian motion (constructed,
using the Brownian web, in \cite{NRS08}) 
and the occurrence of exceptional times for this object.

\subsection{Brownian Web and $(1,2)$ points}
Under diffusive scaling, individual
random walk paths converge to Brownian motions. In \cite{FINR04}, it was proved
(extending the results of~\cite{A81,TW98}) that the entire
collection of discrete paths in the DW converges (in an appropriate
sense) to the continuum Brownian web (BW), which can be loosely
described as the 
collection of graphs of coalescing
one-dimensional Brownian 
motions starting from every possible location in $\R^2$ (space-time).

Formally, the Brownian web (BW) is a random
collection of paths with specified starting points in space-time.
The paths are continuous graphs in a space-time metric space $(\br^2,\rho)$
which is a compactification of $\r^2$.  $(\o,d)$ denotes the space
whose elements are paths with specific starting points. The metric
$d$ is defined as the
maximum of the sup norm of the distance between
two paths and the distance between their respective starting points.
(Roughly, the distance between two paths is small 
when they start from close (space-time)
points and remain close afterwards).
The Brownian web takes
values in a metric space $(\h,d_\h)$,
whose elements are compact collections of paths in
$(\o,d)$ with $d_\h$ the induced Hausdorff metric. Thus the
Brownian web is an
$(\h,\f_\h)$-valued random variable,
where $\f_\h$ is the Borel $\s$-field  associated to the metric
$d_\h$. The next theorem, taken from~\cite{FINR04}, gives some of the key
properties of the BW.
\bteo
\label{teo:char}
There is an \( ({\cal H},{\cal F}_{{\cal H}}) \)-valued random variable
\(
{\W} \)
whose distribution is uniquely determined by the following three
properties.
\begin{itemize}
\item[(o)]  from any deterministic point \( (x,t) \) in
$\r^{2}$,
there is almost surely a unique path \( {B}_{(x,t)} \)
starting
from \( (x,t) \).

\item[(i)]  for any deterministic, dense countable subset
\(
{\cal
D} \) of \( \r^{2} \), almost surely, \( {\W} \) is the
closure in
\( ({\cal H}, d_{{\cal H}}) \) of \( \{ {B}_{(x,t)}: (x,t)\in
{\cal D} \}. \)
\item[(ii)]  for any deterministic  $n$ and \((x_{1}, t_{1}),
\ldots,
(x_{n}, t_{n}) \), the joint distribution of \(
{B}_{(x_{1},t_{1})}, \ldots, {B}_{(x_{n},t_{n})} \) is that
of coalescing Brownian motions (with 
zero drift and unit diffusion
constant).

\end{itemize}
\eteo

This characterization provides a practical construction of the Brownian web. 
For $\cal D$ as above,  construct coalescing Brownian motion paths starting 
from $\cal D$. This defines a {\it skeleton} for the Brownian web. $\W$
is simply defined as the closure of this precompact set of paths. 

\begin{figure}[!ht]
\begin{center}
\includegraphics[width=3cm]{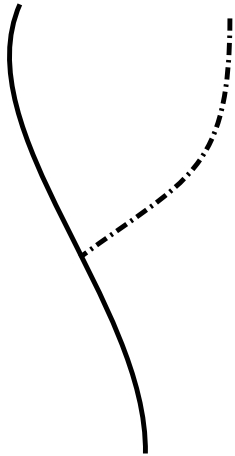}
\caption{A schematic diagram of a left $(m_{\mbox{{\scriptsize\emph{in}}}},
m_{\mbox{{\scriptsize
\emph{out}}}})=(1,2)$ point. 
In this example the incoming path connects to the leftmost
outgoing path,
the right 
outgoing path is a newly born path.
}\label{fig12}
\label{12} 
\end{center}
\end{figure}

We note that generic (e.g., deterministic)
space-time points have almost surely only
$m_{\mbox{{\scriptsize out}}}=1$
outgoing (to later times)
paths from that point
and $m_{\mbox{{\scriptsize in}}}=0$
incoming paths passing through that point (from earlier times). 
An interesting 
property of the BW is related to the 
existence of special points with other values of
$(m_{\mbox{{\scriptsize in}}},m_{\mbox{{\scriptsize out}}})$. 
In the following, a dominant role is played by the $(1,2)$ points as
we shall explain. Back in the lattice,
$(1,2)$ points correspond to 
locations where a path starts
at a ``microscopic'' distance from an old path (that started from earlier time;
we note that in the count of paths, incoming paths that coalesce at
some earlier time are identified)
and coalesces
with it only after some ``macroscopic'' amount of time.
For such a point, the single incident path
continues along exactly
one
of the two outward paths.
The $(1,2)$ point is either left-handed or
right-handed according to
whether the incoming
path connects to the left or right outgoing path.
See Figure~\ref{fig12} for a schematic
diagram of the ``left-handed'' case.
Both varieties occur and it is
known ~\cite{FINR04} that each of the two varieties, 
as a subset of $\R^2$, has Hausdorff dimension~$1$. 

\subsection{The Dynamical Brownian Web and Exceptional Times}
It is natural that there
should also exist scaling limits of the DyDW (including of
the random walk from the origin
evolving in $\tau$, i.e., a dynamical Brownian motion). 
Indeed, this was
proposed by Howitt and Warren~\cite{HW07} 
who also studied (two dynamical time distributional)
properties of any such limit.  
In~\cite{NRS08}, we provided a complete construction
that we now briefly describe.

A priori, a direct construction in the continuum appears difficult since
the DyDW is entirely
based on a modification
of the discrete arrow structure of the DW,
while in the BW it was unclear a priori whether 
there even is
any arrow structure to modify.
Two of the main themes of \cite{NRS08} 
are thus: (i)``Where is the
arrow structure of the BW?'' and (ii)``How is it modified
to yield the DyBW (including a dynamical Brownian motion
from the origin)?''. The answer to
the first question is that the arrow structure of the BW comes
from the $(1,2)$ points. Indeed,
one can change the direction of the ``continuum'' 
arrow at a given $(1,2)$
point $z$ by simply
connecting the incoming path to the newly born path starting from $z$
rather than to the original continuing path. 
(Back in the lattice, this amounts to changing the direction of an arrow
whose switch induces a ``macroscopic'' effect in the web.) 
The answer to question (ii) is based
on the construction of a Poissonian marking of the
$(1,2)$ points (see \cite{NRS08} for details) that indicates which
$(1,2)$ points get switched and at what value of $\tau$
does the switch occur. We note that the main difficulty
in the construction of the DyBW lies in the fact that 
between two dynamical times $\tau<\tau'$,
one needs to switch the direction of a set of $(1,2)$ points
dense in $\r^2$ in order
to deduce the
web at time $\tau'$ from the one at time $\tau$.

\bigskip

We proceed to discuss the existence of exceptional times for $B_0^\tau$,
the dynamical Brownian motion starting from the origin
at dynamical time $\tau$. (We remark that our tameness results,
Theorem~\ref{th1} and Remark~\ref{grw}, can be extended to
the continuum DyBW, but the arguments involve some extra
Brownian web technology.) Recall that the key ingredient
for proving our existence results for the dynamical discrete web is 
contained in Proposition~\ref{sensitivityy-noise}
where we estimate how fast the dynamical discrete web decorrelates. 
The proof of that proposition
mostly relies on the observation that 
$(S_0^\tau,S_0^{\tau'})$ form a sticky pair of random walks. More
precisely, 
we showed in Lemma \ref{pai-sbm} that 
along the $t$-axis the pair alternates between periods during which
the two paths evolve as a single path (they stick) and periods 
during which they move independently.

In \cite{NRS08}, we proved that $\tau\rightsquigarrow B_0^\tau$ 
has a similar structure (as suggested in~\cite{HW07}), 
in that for two distinct dynamical times $\tau,\tau'$, the paths
$B_0^\tau,B_0^{\tau'}$ form 
a $1/(2|\tau-\tau'|)$-sticky pair of Brownian motion. Such a pair can be simply 
expressed in terms of three independent 
standard Brownian motions $(B_d^\tau,B_d^{\tau'},B_s)$ in the following way.
\begin{eqnarray}
\label{dec}
B^\tau_0(t)=B_d^{\tau}(C(t))+B_s(t-C(t)),  \nn\\
B^{\tau'}_0(t)=B_d^{\tau'}(C(t))+B_s(t-C(t)),
\end{eqnarray}
where $C$ is the continuous inverse of
the function
\beqn
\label{dec2}
C^{-1}(s) \, = \, s + \fr{1}{\sqrt{2} \, |\tau-\tau'|}\ l_0(s)
\eeqn
and $l_0$ is the local time at the origin
of the process $\left(B_d^\tau-B_d^{\tau'}\right)/\sqrt{2}$. 
We note that the paths  $B_0^\tau,B_0^{\tau'}$ always spend a strictly positive
Lebesgue measure of time together, hence the name sticky Brownian motions.  
Finally, the time the two paths spend together 
is directly related
to
the parameter
$1/(2|\tau-\tau'|)$ commonly referred to
as the ``amount of stick'' of the pair.

If we denote by $\tilde \pi(\cdot)=\pi(\cdot/\d^2) \d$,
the scaling invariance for the Brownian motion combined with  (\ref{dec}) 
implies that 
$(\tilde B^\tau_0,\tilde B^{\tau'}_0)$ is identical in law to 
a $\d/(\sqrt{2} |\tau-\tau'|)$-sticky pair of Brownian motions. In other words,
the amount 
of stick of the pair $(\tilde B_0^\tau,\tilde B_0^{\tau'})$
vanishes as $\d\rightarrow0$ and from (\ref{dec}) and (\ref{dec2}) we see that
for small $\d$,
\beq
B_0^\tau(t)\approx B_d^{\tau}(t) \ \ \textrm{and} 
\ \ \ B_0^{\tau'}(t)\approx B_d^{\tau}(t),
\eeq
i.e., the two paths become ``almost independent''.
This can be made more precise  by establishing (along the same lines as the proof
of Proposition \ref{sensitivityy-noise}) that
for $$O=\{\forall t\in[0,1], \ \pi(t)>-1 \ \ \textrm{and} 
\ \ \pi(1)>1\}$$ and $\d>0$,
there exist $K,a\in(0,\infty)$ (independent of $\d,\tau$ and $\tau'$) such that 
\beqn
\P(\tilde B^\tau_0\in O \ , \ \tilde B^{\tau'}_0\in O )& \leq & 
\P(\tilde B\in O)^2 \ + \ K (\fr{\d}{|\tau-\tau'|})^a \\
& = & \P(B\in O)^2 \ + \ K (\fr{\d}{|\tau-\tau'|})^a,
\eeqn
where $B$ is a standard Brownian motion.

Since all the results of Section \ref{existence} and Subsection \ref{lowerbound}
for our dynamical random walk
are based on the discrete analog of this result, 
Propositions \ref{violation} and \ref{lower} can easily be extended to
the continuum in the following manner. 
If we define
$$A=\{\inf_{t\in[0,1]}B(t)>-1\ , \ B(1)>1\}$$
and let
$\g_0=1/\P(A)$, we have:

\bteo
Let $\bar {\mathcal T}(K)$ be the set
of $\tau$'s belonging to $[0,\infty)$ such that
\beq
\forall t\geq0 , \ \ B^\tau_0(t)\geq -1 -K\sqrt{t}.
\eeq
Then $\bar {\mathcal T}(K)$ is non-empty and
\begin{equation}
dim_H(\bar {\mathcal T}(K)) \ \geq \ 1\, - \,
\frac{\log{\gamma_0}}{\log{\gamma(K)}}
\ \textrm{for} \ K\, > \, K(\g_0) \, .
\end{equation}
Thus, $\lim_{K\uparrow\infty} dim_H(\bar {\mathcal T}(K))=1$.
\eteo

We conclude by noting that also our upper bound results on
the Hausdorff dimension, Propositions~\ref{upper} and~\ref{tme-K},
can be extended to the continuum DyBW, 
but, like the tameness results,
that extension requires some extra Brownian web technology
beyond what 
is described in this paper. 

\appendix
\section{Some Estimates On Random Walks\label{appendix} 
(Proof of Lemma \ref{principal22})}

We start with the two following lemmas.
\blem\label{S777}\cite{S77}
\label{samere}
Let $j,K\in(0,\infty)$ 
and let $B$ be a standard Brownian motion. Then there exists
$q \in (0,\infty)$ such that
\begin{equation}
\lim_{t \to \infty}  \ \ {t}^{p(K)/2} \mathbb{P}(\forall s\in [0,t], 
\ B(s)\geq-j-K\sqrt{s}) = \ q \, ,
\end{equation}
where $p({K})$ is the solution in $(0,1)$ of
the equation
\begin{equation}
f(p,{K}) \equiv
\frac{\sin(\pi p / 2) \Gamma(1+ p / 2)}{\pi} \sum_{n=1}^{\infty}
\frac{(\sqrt{2} K)^n}{n !} \ \Gamma((n-p)/2)=1 \, .
\end{equation}
Furthermore, $p(K)$ is a continuous decreasing function
on $(0,\infty)$ with
\begin{equation}
\lim_{K\uparrow\infty}p(K)=0 \ \ \ \txt{and} \ \
\lim_{K\downarrow0}p(K)=1.
\end{equation}
\elem

\blem
\label{enfinlafin}
Let $K\in(0,\infty)$, $ l\in(0,1)$ and let $S$ be a simple symmetric
random walk. Then there exists
$\bar c(K,l) \in (0,\infty)$ such that for every $n$
\beq
{n}^{p(K/l)/2} \ \P(\forall k\leq n, \  \ S(k)\geq-1-K\sqrt{k})\  \leq \  \bar c(K,l).
\eeq
\elem
\begin{proof}
By Lemma \ref{S777}, it suffices to prove that for every 
$l<1$, there exists $c(K,l)$ such that
\beq\label{p1,p2}
\P(\forall k\leq n,\ S(k)\geq -1-K\sqrt{k}) \leq 
\   c(K,l) \ \  \P(\forall t\in [0,l^2n],\ B(t)\geq -2-K\sqrt{t}/l).
\eeq
We now prove the latter inequality.
Consider $S$ the discrete time random walk embedded in
the Brownian motion $B$. Namely,
we define inductively a sequence of stopping times
$t_i$  with $t_{0}=0$ and
\begin{equation}
\label{tautau}
t_{i+1}=\inf\{t>t_i :  \
|B(t)- B(t_i)| \geq 1\}
\end{equation}
and then we define
$S(i)=B(t_i)$. Note that $S$ and $\{t_i\}$ are independent and therefore
\begin{eqnarray*}
&&\P(\forall k\leq n, \ S(k)\geq-1-K\sqrt{k})  \\[0.1in] 
&&=\P(\forall k\leq n, \ S(k)\geq-1-K\sqrt{k} , \ l^2 k  
\leq \ t_k\leq \fr{k}{l^2}) \ / \ \P(\forall k\leq n,\ l^2 k \leq 
\ t_k\leq \fr{k}{l^2} ), \\
&&=  \P(\forall k\leq n, \ B(t_k)\geq-1-K\sqrt{k} , \ l^2 k  
\leq \ t_k\leq \fr{k}{l^2}) \ / \ \P(\forall k\leq n,\ l^2 k  \leq 
\ t_k\leq \fr{k}{l^2} ), \\
&& \leq \P(\forall k\leq n, \ B(t_k)\geq-1-K\sqrt{t_k}/l , \ l^2 k  
\leq \ t_k) \ / \ \P(\forall k\leq n,\ l^2 k  \leq \ t_k\leq \fr{k}{l^2}).
\end{eqnarray*}
If for every $k\leq n$, we have $B(t_k)\geq-1-K\sqrt{t_k}/l$ and moreover 
$l^2 k \leq t_k$, then on $[0,l^2 n]$, everytime $B$ takes an integer 
value, $B$ is to the right of $t\rightsquigarrow-1-K\sqrt{t}/l$. Hence, 
$B$ remains to the right of $t\rightsquigarrow-2-K\sqrt{t}/l$ on $[0,l^2 n]$ 
which implies that
\begin{eqnarray*}
&& \P(k\leq n, \ S(k)\geq-1-K\sqrt{k})\\ 
&& \leq 
\P(\forall t\leq l^2 n, \  B(t)\geq-2-K\sqrt{t}/l) \ / \ 
\P(\forall k\leq n,\ l^2 k  \leq \ t_k\leq \fr{k}{l^2}).
\end{eqnarray*}
Finally, $t_k$ is a sum of $k$ i.i.d. random variables with mean $1$
(whose common distribution includes $1$ in its support). Therefore, 
$\P(\forall k\in\N,\ l^2 k  \leq \ t_k\leq \fr{k}{l^2})>0$ 
and (\ref{p1,p2}) follows. 
This completes the proof of the lemma.

\end{proof}

We are now ready to prove Lemma \ref{principal22}.
Recall that $S_\e$ is a simple
random walk with
$$\P(S_\e(k+1)-S_{\e}(k)=1)=\fr{1}{2}+\fr{1}{2}(1-e^{-\e}).$$
Since $\fr{1}{2}+\fr{1}{2}(1-e^{-\e})\leq \fr{1}{2}(1+\e)$, 
it is enough to show the conclusions of the lemma for the simple walk 
$\bar S_\e$ where $p_\e^{\pm}=\P(\bar S_\e(k+1)-\bar S_{\e}(k)=\pm1)=
\fr{1}{2}(1\pm\e).$

Let $T_\e=\inf\{n>0 :\bar S_\e(n)<-1- K\sqrt{n}\}$. We have 
\beqn
\P(T_\e=n)
& = & \P(T_0=n) \  f_\e(n) \\[0.1in]
\txt{with} & &  f_\e(n)\equiv 
({2 p_\e^-})^{\fr{1}{2}(n+\lfloor 1+K\sqrt{n} \rfloor +1)} \ 
({2 p_\e^+})^{\fr{1}{2}(n-\lfloor 1+K\sqrt{n}\rfloor -1)}, \eeqn
where $\lfloor x \rfloor$ denotes the greatest integer $\leq x$.
Since a simple symmetric random walk $S$ eventually
hits the moving boundary $t \rightsquigarrow-1-K\sqrt{t}$ we have
\beqn
\P(T_\e=\infty)=1-\sum_{n\geq 1} f_\e(n) \P(T_0=n) \ = 
\ \sum_{n\geq 1}(1-f_\e(n)) \P(T_0=n).
\eeqn
Then, proceeding to a summation by parts we have
\beqn
\P(T_\e=\infty)
&=&
\sum_{n\geq 1} (\P(T_0\geq n)-\P(T_0\geq n+1)) \ [1-f_{\e}(n)], \\[0.1in]
&=&
\sum_{n\geq 1} \P(T_0\geq n+1)(f_\e(n)-f_\e(n+1)) 
\ + \ (1-f_{\e}(1)),\\[0.1in]
&= & \S_\e \ + \ (1-f_{\e}(1)). \label{A102}\\[0.1in]
\txt{with} & & \S_\e \equiv \sum_{n\geq 1} 
\P(T_0\geq n+1)f_\e(n)(1-\fr{f_\e(n+1)}{f_\e(n)})
\eeqn
We proceed to estimate $\S_\e$. First, for $\e\in(0,1)$, we have
\begin{eqnarray*}
\ln(1+\e)\leq {\e}, \ \ \ \ln(1-\e)\leq -{\e},
\end{eqnarray*}
implying that
\begin{eqnarray*}
f_\e(n)& = &\exp\{\fr{1}{2}\ln(1-\e)(n+ \lfloor 1+K\sqrt{n} \rfloor +1)\ + 
\ \fr{1}{2}\ln(1+\e)(n- \lfloor 1+K\sqrt{n} \rfloor -1)\} \\[0.1in]
&\leq & 
\exp\{-\fr{\e}{2}(n+ \lfloor 1+K\sqrt{n} \rfloor +1) \ + 
\ \fr{\e}{2}(n- \lfloor 1+K\sqrt{n} \rfloor -1)\} \\[0.1in]
&\leq & \exp\{ -\e ( \lfloor 1+K\sqrt{n} \rfloor +1) \} \\[0.1in]
&\leq & \exp\{-{\e} K\sqrt{n}\}.
\end{eqnarray*}
Next, if we set 
$\D_n\equiv \lfloor 1+K\sqrt{n+1} \rfloor - \lfloor 1+K\sqrt{n} \rfloor 
= \lfloor K\sqrt{n+1} \rfloor - \lfloor K\sqrt{n} \rfloor $, we have
\begin{eqnarray*}
1-\fr{f_\e(n+1)}{f_\e(n)}& = &1 \ - \  \exp\{\ln(1-\e)(\fr{1}{2}+
\fr{1}{2}\D_n)+\ln(1+\e)(\fr{1}{2}-\fr{1}{2}\D_n)\}, \nonumber\\
& = &1 \ - \ \exp(-\e \D_n-\fr{\e^2}{2} + o(\e) \D_n+ o(\e^2)) = 
\e \D_n +\fr{\e^2}{2} + o(\e) \D_n+ o(\e^2) .
\end{eqnarray*}

By Lemma \ref{enfinlafin}, for every $l< 1$ there exists $\bar c(K,l)$ such that
\begin{eqnarray}
\S_\e & \leq \sum_{n\geq 1} \P(T_0\geq n+1) \ \exp\{-{\e} K\sqrt{n}\} \ 
(\e \D_n +\fr{\e^2}{2} + o(\e) \D_n \ + \       o(\e^2)) \\[0.15in]
& \leq \bar c(K,l) \sum_{n\geq 1} \exp\{-{\e} K\sqrt{n}\} \ 
(\e \fr{\D_n}{{n}^{p/2}} +\fr{\e^2}{2{n}^{p/2}} + 
o(\e) \fr{\D_n}{{n}^{p/2}} + \fr{o(\e^2)}{{n}^{p/2}}) \label{enfinlafin3}.
\end{eqnarray}
where $p\equiv p(K/l)$ is as in Lemma \ref{enfinlafin}. Since 
$\D_n= \lfloor K\sqrt{n+1} \rfloor - \lfloor K\sqrt{n} \rfloor $, 
and $(\sqrt{n+1}-\sqrt{n})\sqrt{n}\raw 1/2$ 
as $n\raw\infty$, 
it is natural to expect that
\beqn
&&\lim_{\e \daw0} \ \e^2 \ \sum_{n\geq 1} \ \exp\{-{K}\sqrt{\e^2 n}\} \ \fr{\D_n}{\e} \ \fr
{1}{(\sqrt{\e^2 n})^p} \ \nn  \\[0.1in]
& = & 
\lim_{\e \daw0} \ \e^2 \ \sum_{n\geq 1} \ \exp\{-{K}\sqrt{\e^2 n}\} \ 
\fr{K}{2\sqrt{n \e^2}} \ \fr{1}{(\sqrt{\e^2 n})^p} \nonumber \\
& = &\ \fr{1}{2} \ \int_{0}^\infty  {\exp\{-{K \sqrt{t}}}\}\ 
\fr{K}{ {t}^{{(p+1)}/2}} \ dt\label{enfinlafin1}
\eeqn
where the second equality is due to the Riemann sum on the right hand side 
of the first equality.
To justify the first equality one may note that $\D_n$ is (for large $n$)
either $0$ or $1$ and then define $N_\ell (n)$ (resp., $N_u (n)$) to
be the largest $m\leq n$ (resp., smallest $m>n$) such that $\D_m \neq 0$.
It is straightforward to show first that that $N_u (n)-N_\ell (n)/\sqrt{n} \to 2/K$
as $n \to \infty$ and then to obtain the first equality of~(\ref{enfinlafin1}) as
a consequence.
It is also the case that
\beq
\lim_{\e \daw0} \ \e^2 \ \ \sum_{n\geq 0} \  \exp\{-{K}\sqrt{\e^2 n}\} \ 
\fr{1}{{(\e^2 n)}^{p/2}} = \  \ \int_{0}^\infty  \exp\{-K \sqrt{t}\}
\fr{1}{ {t}^{p/2}} \ dt \label{enfinlafin2}.
\eeq
Since $0<p<1$ both integrals in (\ref{enfinlafin1}) 
and (\ref{enfinlafin2}) are finite. Thus,
(\ref{enfinlafin3})
yields $S_\e=O(\e^p)=O(\e^{p(K/l)})$.

Finally, it is easy to prove that $f_\e(1)-1=O(\e)$. Since $p(K/l)<1$, 
Lemma \ref{principal22} follows from
(\ref{A102}).

\bigskip
\bigskip

{\em Acknowledgements.}
The research of L.R.G.~Fontes was supported in part by FAPESP grant
2004/07276-2 and CNPq grants 307978/2004-4 and 484351/2006-0;
the research 
of the other authors was supported in part by N.S.F. grants DMS-01-04278 and 
DMS-06-06696.

\end{document}